\documentclass{amsart}
\usepackage{amsmath,amssymb}
\usepackage[dvips]{graphicx} 
\usepackage[arrow,matrix]{xy}

\newtheorem{thm}{Theorem}[section]    
   
\newtheorem{prop}[thm]{Proposition}

\newtheorem{conj}[thm]{Conjecture}

\theoremstyle{definition}
\newtheorem{defn}[thm]{Definition}
\newtheorem{example}[thm]{Example}

\theoremstyle{remark}
\newtheorem{rem}[thm]{Remark}

\numberwithin{equation}{section}

\newcommand\bC{\mathbb C}
\newcommand\bZ{\mathbb Z}

\newcommand\bP{\mathbb P}
\newcommand\bQ{\mathbb Q}

\newcommand\cV{\mathcal V}
\newcommand\cL{\mathcal L}
\newcommand\cT{\mathcal T}

\newcommand\cD{\mathcal D}
\newcommand\cO{\mathcal O}
\newcommand\wcD{\widehat{\mathcal D}}
\newcommand\Ztwo{\bZ/2}

\newcommand{\Ch}{\operatorname{Ch}}
\newcommand{\Ind}{\operatorname{Ind}}

\newcommand{\sign}{\operatorname{sign}}

\newcommand{\Ca}{$C^*$-algebra}
\newcommand\lp{\textup{(}}
\newcommand\rp{\textup{)}}
\newcommand\wA{\widehat A}
\newcommand\co{\colon\,}
\newcommand\wV{\widehat V}
\newcommand\wZ{\widehat Z}
\newcommand\wM{\widehat M}
\newcommand\wi{\widehat \iota}

\begin{document}
\title[A Novikov Conjecture in Algebraic Geometry]{An analogue
of the Novikov Conjecture\\%
in complex algebraic geometry}                    
\author{Jonathan Rosenberg} 
\address{Department of Mathematics\\University of Maryland\\College
Park, MD 20742, USA}
\email{jmr@math.umd.edu}
\urladdr{http://www.math.umd.edu/\lower.6ex\hbox{\char"7E}jmr}
\thanks{This work was partially supported by NSF grant number
DMS-0504212.} 

\date{}

\dedicatory{}


\begin{abstract}
We introduce an analogue of the Novikov Conjecture on higher
signatures in the context of the algebraic geometry of (nonsingular)
complex projective varieties. This conjecture asserts that
certain ``higher Todd genera'' are birational invariants. This implies
birational invariance of certain extra combinations of Chern classes
(beyond just the classical Todd genus)
in the case of varieties with large fundamental group (in the
topological sense). We prove the conjecture under the assumption of
the ``strong Novikov Conjecture'' for the fundamental group, which is
known to be correct for many groups of geometric interest.
We also show that, in a certain sense, our conjecture is best possible.
\end{abstract}
\maketitle

%
%

\section{Introduction and Statement of the Conjecture}
\label{sec:intro}

In recent years, there has been considerable interest in the Novikov
Conjecture on higher signatures, in its variants and analogues, and on
methods of proof coming from index theory and noncommutative
geometry. (See for example the books \cite{FRR1}, \cite{FRR2}, and
\cite{KrL}, and 
papers such as \cite{K1} and \cite{KM}. More references to the vast
literature may be found in the books just cited.) The classical
statement of the Novikov Conjecture is as follows:
\begin{conj}[Novikov]
\label{conj:NC}
Let $M^n$ be a closed oriented manifold, with fundamental
class $[M]\in H_n(M,\bZ)$ and total $L$-class $\cL(M)\in H^*(M,\bQ)$. Let
$\pi$ be a discrete group, let $x\in H^*(B\pi, \bQ)$,
and let $u\co M\to B\pi$ be a map from $M$
to the classifying space of this group. {\lp}In practice, usually
$\pi=\pi_1(M)$ and $u$ is the classifying map for the universal cover
of $M$, though it is not necessary that this be the case.{\rp} Then
the ``higher signature'' $\langle \cL(M)\cup u^*(x),\,[M]\rangle$ is
an oriented homotopy invariant of $M$. In other words, if $\wM$ is
another closed oriented manifold, and if $h\co \wM\to M$ is an
orientation-preserving homotopy equivalence, then
\[
\langle \cL(M)\cup u^*(x),\,[M]\rangle = \langle \cL(\wM)\cup
h^*\circ u^*(x),\,[\wM]\rangle .
\]
\end{conj}
Since $\cL$ is a power series in the Pontrjagin classes, this
conjecture implies an additional rigidity for the Pontrjagin classes
of manifolds with ``large'' fundamental group, beyond the constraint
of homotopy-invariance of the Hirzebruch signature $\sign M = \langle
\cL(M), [M]\rangle $. For example, if $M$ is aspherical, the
conjecture implies that \emph{all} of the rational Pontrjagin classes
of $M$ are homotopy invariants, which is certainly not the case for
simply connected manifolds, for example.

We shall in this paper consider an analogue of the Novikov Conjecture
in complex algebraic geometry, with higher signatures replaced by
higher Todd genera, and with homotopy invariance replaced by
birational invariance. I would like to thank Professor Nigel Hitchin
of Oxford University for a conversation (in 1998) that led
to the following statement:

\begin{conj}
\label{conj:higherTodd}
Let $V$ be a nonsingular complex projective variety of complex
dimension $n$, also viewed as a complex manifold.
Let $[V]\in H_{2n}(V,\bZ)$ be its fundamental
class and $\cT(M)\in H^*(V,\bQ)$ be its total Todd class
in the sense of \textup{\cite[Ch.\ 3]{Hir}}. {\lp}Homology and cohomology
are taken in the ordinary sense, for the underlying complex
manifold of $V$ with its locally compact Hausdorff topology.{\rp}
Let $\pi$ be a discrete group, let $x\in H^*(B\pi, \bQ)$,
and let $u\co V\to B\pi$ be a map from $V$
to the classifying space of this group.  Then
the ``higher Todd genus'' $\langle \cT(V)\cup u^*(x),\,[V]\rangle$ is
a birational invariant of $V$ {\lp}within the class of
nonsingular complex projective varieties{\rp}.
\end{conj}

Clearly this statement is very closely analogous to Conjecture
\ref{conj:NC}. It seems appropriate to use birational maps in this
context, as they provide the appropriate notion of weak isomorphism in
the category of projective varieties.
Conjecture \ref{conj:higherTodd} again implies an extra rigidity for the
Chern classes of $V$ in the case where the fundamental group of $V$
is ``large,'' beyond the birational invariance of the Todd
genus $\langle \cT(V), [V]\rangle$, which is equal to the arithmetic
genus $\sum_{j=1}^n (-1)^j \dim H^j(V, \cO_V)$ \cite[Theorem 20.2.2]{Hir}.
(Here and throughout the paper I mean fundamental group
in  the usual topologists' sense, not in the sense of algebraic geometry.)

Most of the literature on the Novikov Conjecture (cf.\ \cite{FRR1},
\cite{FRR2}, and \cite{KrL}) does not deal with
the statement of Conjecture \ref{conj:NC} directly. Instead, this
statement is reduced, using either surgery theory or index theory, to
a statement about injectivity of a certain \emph{assembly map}, which
is a property of the group $\pi$ and which does not involve the
manifold $M$ at all.  Different formulations are possible, but here we
will use:
\begin{defn}[\textup{\cite{R}}]
\label{def:strongNC}
Let $\pi$ be a discrete group, and let $C^*(\pi)$ be its group {\Ca}
{\lp}either full or reduced, it doesn't matter{\rp}. Define the
\textbf{assembly map} for $\pi$, $A_\pi\co K_*(B\pi)\to K_*(C^*(\pi))$, as
follows. Let $\cV_{B\pi}$ be the \textbf{universal flat
$C^*(\pi)$-bundle} over $B\pi$, $\cV_{B\pi} = E\pi \times_\pi
C^*(\pi)$. This is a bundle of rank-one projective {\lp}right{\rp}
$C^*(\pi)$-modules over $B\pi$. As such, it has a class $[\cV_{B\pi}]$ in
the Grothendieck group $K^0(B\pi; C^*(\pi))$. (See \cite{MF}
or \cite[\S1.3]{MR1787114}.) Then $A_\pi$ is the 
Kasparov product (a kind of slant
product) with this class. (If there is no compact model for $B\pi$,
we interpret the $K$-homology $K_*(B\pi)$ as $\varinjlim K_*(X)$, as
$X$ runs over the finite subcomplexes in a CW model for $B\pi$, and
take the inductive limit of the slant products with $[\cV_X]\in K^0(X;
C^*(\pi)) \cong K_0(C(X)\otimes C^*(\pi)$.)

Then we say that $\pi$ satisfies the \textbf{Strong Novikov
Conjecture} (SNC) if $A_\pi$ is rationally injective.
\end{defn}
It was shown by Kasparov (\cite{K1}, \cite{K2}, \cite{K3}) --- see
also \cite{KM} for a variant of the proof --- that the Strong Novikov
Conjecture for $\pi$ implies Conjecture \ref{conj:NC} (for the same
$\pi$, and for arbitrary manifolds $M$). The SNC also implies \cite{R} the
non-existence of metrics of positive scalar curvature on closed spin
manifolds for which ``higher $\wA$-genera'' are non-zero.

The Strong Novikov Conjecture is weaker than the Baum-Connes
Conjecture and is now known to hold for many
(overlapping) classes of groups of
geometric interest: discrete subgroups of Lie groups \cite{K1}, groups
which act properly on ``bolic'' spaces \cite{KS}, groups with finite
asymptotic dimension \cite{Yu1}, groups with a uniform embedding into
a Hilbert space \cite{Yu2}, and hyperbolic groups \cite{MinYu}.

The main result of this paper will be
\begin{thm}
\label{thm:main}
The Strong Novikov Conjecture {\lp}of Definition
\ref{def:strongNC}\,{\rp} for a group $\pi$ implies Conjecture
\ref{conj:higherTodd}, for all nonsingular complex projective
varieties and for the same group $\pi$.
\end{thm}

This result is thus formally analogous to Kasparov's.  The proof will
be given in Sections \ref{sec:proofmethod} and
\ref{sec:proofdetails} of this paper, and will involve
mimicking certain aspects of one of the proofs of birational invariance
of the Todd genus.

Thus Conjecture \ref{conj:higherTodd} is true as long as the group
$\pi$ is a discrete subgroup of a Lie group, a group with finite
asymptotic dimension, a group with a uniform embedding into
a Hilbert space, or a hyperbolic group. This implies some extra
birational
rigidity statements for Chern classes of certain smooth projective
varieties; we will give some concrete examples in Section
\ref{sec:concl}. As we shall see in the Section
\ref{sec:remarks}, there is also a sense in
which Conjecture \ref{conj:higherTodd} is best possible. 

After this paper was written and posted in August, 2005, 
J\"org Sch\"urmann kindly pointed out to me
that Conjecture \ref{conj:higherTodd},
even without assuming the strong Novikov Conjecture, can be deduced from
the result in \cite[Example 3.3, part 3]{BSY} that 
when $f\co X \to V$ is a morphism of
smooth projective varieties and a birational equivalence,
then $f_* (\cT(X)\cap [X]) = \cT(V) \cap [V]$. A similar argument
can be found in \cite{BW}, which points out that the results
of \cite{MR549773} can be used to prove the corresponding fact in $K$-homology,
that under the same hypotheses, $f_*([D_X]) = [D_V]$. (The notation
in this equation will be explained in section \ref{sec:proofmethod}.)
This implies the $K$-homology analogue of Conjecture \ref{conj:higherTodd},
that $u_*([D_V])\in K_0(B\pi)$ is a birational invariant of $V$.
This fact is somewhat stronger than Conjecture \ref{conj:higherTodd},
since it implies the latter but also gives some torsion information.
And finally, Borisov and Libgober \cite{BoL} have now proven an analogue
of Conjecture \ref{conj:higherTodd} with the Todd genus replaced by
the elliptic genus (which includes it as a special case) and with
birational invariance replaced by $K$-invariance.  (Two smooth
complex projective varieties are said to be $K$-equivalent --- see
\cite{Wang} for a nice survey --- if they
are birationally equivalent and if they also have ``the same $c_1$,''
more exactly if they have a common blowup such that the pull-backs 
of their canonical line bundles to this common blowup are
equivalent.)

Thus it is fair to say that this paper is now somewhat superseded.
Nevertheless, we have decided to publish it anyway in what is essentially
its original form.  There are a few reasons for this. First of all, the
appearance of \cite{BW} and of  \cite{BoL} confirms that our
Conjecture \ref{conj:higherTodd} seems to be of interest.  Secondly,
the results of section \ref{sec:remarks}
appear to be new and of independent interest.
And finally, the method of proof, and in particular 
Theorem \ref{thm:holAbundle},
may prove to be useful for other applications.

I would like to thank many friends and colleagues for discussions
or correspondence on the subject of this paper. These include
Paulo Aluffi, Jonathan Block, Nigel Higson,
Jack Morava, Niranjan Ramachandran, J\"org Sch\"urmann,
Shmuel Weinberger, and many others.


\section{Method of Proof of the Main Theorem}
\label{sec:proofmethod}

Before we begin the proof of Theorem \ref{thm:main}, we need first to
clarify precisely what it means. If you look carefully, you will
notice that Conjecture \ref{conj:higherTodd} is not quite as precise as
Conjecture \ref{conj:NC}, the reason being that a birational ``map''
is only defined in the complement of a subvariety of lower
dimension. In other words, it
is not necessarily a map (in the sense of topologists) at all, only a
correspondence. However, fortunately we can appeal to:
\begin{thm}[Factorization Theorem \cite{Wod}, \cite{AKMW}] 
\label{thm:factorization}
Any birational map 
\[ \phi\co \xymatrix{\wV\ar@{.>}[r] & V}
\]
between nonsingular complex projective varieties $\wV$ and $V$ can be
factored into a series of blowings up and blowings 
down with nonsingular irreducible centers. Furthermore, if $\phi$ is
an isomorphism on an open subset $U$ of $V$, then the centers of the
blowings up and blowings can be chosen disjoint from $U$, with all the
intermediate varieties projective.  
\end{thm}

As a consequence, to make sense of Conjecture \ref{conj:higherTodd},
it suffices to deal with the case of a ``blowing down'' $\phi\co \wV
\to V$, which is an actual morphism of varieties, and in particular a
smooth map of complex manifolds. We may assume (by Theorem
\ref{thm:factorization}) that $\wV$ is obtained from
blowing up $V$ along a nonsingular irreducible subvariety $Z$ of (complex)
codimension at least $2$. In that event, $\phi$ restricts an isomorphism
$(\wV\smallsetminus \wZ) \xrightarrow{\cong} (V\smallsetminus Z)$,
and $\wZ = \phi^{-1}(Z)$ is a nonsingular subvariety of $\wV$ of
(complex) codimension $2$. Thus we are reduced to showing:

\begin{conj}
\label{conj:precisehigherTodd}
Let $\phi\co \wV \to V$ be a blowing down of nonsingular complex
projective varieties as above, let $\pi$ be a discrete group, let
$x\in H^*(B\pi, \bQ)$, and let $u\co V\to B\pi$ be a map from $V$
to the classifying space of this group.  Then
\[
\langle \cT(V)\cup u^*(x),\,[V]\rangle = \langle \cT(\wV)\cup
\phi^*\circ u^*(x),\,[\wM]\rangle .
\]
\end{conj}
This statement is now precisely analogous to Conjecture \ref{conj:NC}.

Incidentally, it is worth pointing out a fact well known to experts,
though not entirely trivial:
\begin{prop}
\label{prop:IsoOnPi}If $\phi\co \wV \to V$ is a blowing down of
nonsingular complex projective varieties as above, then as a map of
complex manifolds, $\phi$ induces an isomorphism on fundamental
groups. 
\end{prop}
\begin{proof}
Since the open subset $V\smallsetminus Z$ is the complement of $Z$,
which has complex codimension at least $2$
(and thus real codimension at least $4$), a general position
argument shows that the inclusion $\iota\co V\smallsetminus Z
\hookrightarrow V$ induces an isomorphism on $\pi_1$. Similarly, since
$\wZ$ has complex codimension $1$ (and thus real codimension $2$) in
$\wV$, $\wi\co \wV\smallsetminus \wZ
\hookrightarrow \wV$ induces an epimorphism on $\pi_1$, though not
necessarily (just from general position arguments alone) an
isomorphism. Now chase the commutative diagram 
\[
\xymatrix{
\pi_1(\wV\smallsetminus \wZ) \ar[d]_{\phi_*}^{\cong} \ar@{->>}[r]^(.55){\wi_*} 
& \pi_1(\wV) \ar[d]^{\phi_*}\\
\pi_1(V\smallsetminus Z) \ar[r]^(.55){\iota_*}_(.55)\cong & \,\pi_1(V).}
\]
Since $\iota_*\circ \phi_*$ is surjective, the map $\phi_*$ on the
right is surjective, i.e., $\phi_*\co \pi_1(\wV) \to \pi_1(V)$ is an
epimorphism. But since $\iota_*\circ \phi_*$ is also injective and
$\wi_*$ is surjective,
$\phi_*\co \pi_1(\wV) \to \pi_1(V)$   must be an isomorphism.
\end{proof}

Proposition \ref{prop:IsoOnPi} at least makes it plausible that there
should be birational invariants of projective varieties coming from
the cohomology of the fundamental group.  In fact, this Proposition
proves Conjecture \ref{conj:precisehigherTodd}, and thus
Conjecture \ref{conj:higherTodd},  in the special case
where $u\in H^{2n}(B\pi, \bQ)$, $n$ the (complex) dimension of $\wV$
and $V$, since the constant term in both $\cT(V)$ and $\cT(\wV)$ is
$1$. 

Now we can explain the strategy of the proof of Theorem
\ref{thm:main}. Suppose $\phi\co \wV \to V$ is a blowing down of
nonsingular complex projective varieties as in Proposition
\ref{prop:IsoOnPi}, and fix a map $u\co V \to B\pi$. Then we can use
$u$ to pull back the universal flat
$C^*(\pi)$-bundle over $B\pi$, $\cV_{B\pi}$ (see Definition
\ref{def:strongNC}), to a flat $C^*(\pi)$-bundle $\cV_V$ over $V$. Fix
K{\"a}hler metrics on $V$ and $\wV$, and let $D_V$ and $D_{\wV}$ be
the associated Dolbeault operators $\bar\partial + {\bar\partial}^*$,
acting on forms of type $(0,*)$. These operators are odd with respect
to the $\bZ/2$-grading of the form bundle by parity of degree.
A flat connection on $\cV_{B\pi}$
pulls back to flat connections on $\cV_V$ and on $\phi^*(\cV_V) =
\cV_{\wV}$. Using these connections, we can make sense of ``$D_V$ with
coefficients in $\cV_V$,'' say $\cD$, and of ``$D_{\wV}$ with
coefficients in $\cV_{\wV}$,'' say $\wcD$. These are ``elliptic operators with
coefficients in a {\Ca}'' in the sense of Mishchenko and Fomenko
\cite{MF,MR1787114}. Again, they are odd with respect
to the $\bZ/2$-grading of the bundles on which they act.
As such, they have ``indices'' in the sense of \cite{MF}
which take values in the $K$-group
$K_0(C^*(\pi))$. (The rough idea is that each of these
operators should have a kernel which is a $\Ztwo$-graded finitely
generated projective module over $C^*(\pi)$, and the ``index'' is the
formal difference of the classes of the even and odd
parts of the kernel in the Grothendieck
group of finitely generated projective modules, i.e., in
$K_0(C^*(\pi))$. Strictly speaking, things are slightly more
complicated, since the kernel of an elliptic operator in the sense of
Mishchenko and Fomenko need not be projective. But the statement above is true
after ``compact perturbation'' of the operator, and one can show
that the choice of compact perturbation does not affect the index.
See \cite{MF} and \cite[\S1.3.4]{MR1787114} for more details. )

Next, we apply the index theorem of Mishchenko and
Fomenko \cite{MF} to this situation. Exactly as in \cite{R} (the only
difference being that here we are using the Dolbeault operator; there
we were using the signature operator), the result is that
\begin{multline}
\label{eq:indexcalc}
\Ind(\cD) = \langle [D_V], [\cV_V] \rangle =
\langle [D_V], u^* [\cV_{B\pi}] \rangle \\ = \langle u_*[D_V],
[\cV_{B\pi}] \rangle = A_\pi([V \xrightarrow{u} B\pi])
\in K_0(C^*(\pi)),
\end{multline}
where $A_\pi$ is the assembly map. Here $[D_V]\in K_0(V)$ is the
$K$-homology class 
of the Dolbeault operator $D_V$ in the sense of Kasparov (\cite{K0}
and \cite{K}), $[\cV_V]$ is the class of the $C^*$-bundle $[\cV_V]$
in the Grothendieck group $K^0(V; C^*(\pi))\cong K_0(C(V)\otimes
C^*(\pi))$ of such bundles, $\langle\phantom{X},\,\phantom{X}\rangle$
is the pairing between $K$-homology and $K$-cohomology, and $[V
\xrightarrow{u} B\pi]\in K_0(B\pi)$ is the 
$K$-homology class defined by the map $u$ from the complex manifold
$V$ into $B\pi$. (This $K$-homology class is the image of the class of $V
\xrightarrow{u} B\pi$ in complex bordism $\Omega^U_{2n}(B\pi)$ under the
map $\Omega^U_{2n}(B\pi)\to K_{2n}(B\pi)\cong K_0(B\pi)$ given by
the map of spectra $\mathbf{MU} \to \mathbf{K}$ defined by the
usual $K$-theory orientation of complex manifolds---see
\cite[Ch.\ VII]{Stong}.)

If we now rationalize, applying the Chern character $\Ch$ (which is an
isomorphism from rationalized $K$-cohomology and $K$-homology to
ordinary rational cohomology and homology) to equation
\eqref{eq:indexcalc} and recalling the
fact that $\Ch([D_V])$ is Poincar\'e dual to the Todd class
(see for example \cite[\S4]{AS3}), we obtain:
\begin{multline}
\label{eq:rationalindexcalc}
\Ind(\cD) \otimes \bQ = \langle \cT(V)\cup \Ch[\cV_V], [M]\rangle
= \langle \cT(V)\cup u^* \Ch [\cV_{B\pi}], [M]\rangle \\
= \langle \cT(V)\cap [M],  u^* \Ch [\cV_{B\pi}] \rangle
\in K_0(C^*(\pi))\otimes \bQ.
\end{multline}
The Strong Novikov Conjecture implies that all the 
even-dimensional rational
cohomology of $\pi$ can be obtained from the class 
$\Ch [\cV_{B\pi}]\in H^*(B\pi, K_0(C^*(\pi))\otimes \bQ)$
(via homomorphisms $K_0(C^*(\pi))\otimes \bQ \to \bQ$),
and thus that all higher Todd genera (in the sense of Conjecture
\ref{conj:higherTodd}) can be obtained from $\Ind(\cD)$.
So to prove Theorem \ref{thm:main}, it will suffice to prove the
following slightly stronger statement:
\begin{thm}
\label{thm:indexinvariance}
Let $\phi\co \wV \to V$ be a blowing down of nonsingular 
$n$-dimensional complex
projective varieties as above, and let $\pi$ be a discrete group.
Let $\cD$ and $\wcD$ be the $C^*(\pi)$-linear Dolbeault operators on
$V$ and on $\wV$ as above. Then $\Ind(\cD)=\Ind(\wcD)$ in
$K_0(C^*(\pi))$. 
\end{thm}

We will prove Theorem \ref{thm:indexinvariance} by reducing it to
equality of certain spaces of holomorphic forms on $\wV$ and $V$,
using the fact that these manifolds agree in the complement of the
lower-dimensional varieties $\wZ$ and $Z$.

\section{Technical Details}
\label{sec:proofdetails}

The proof of Theorem \ref{thm:indexinvariance} proceeds in several
steps. Some of these steps have to do with technicalities of
Mish\-chenko-Fo\-men\-ko index theory; others would be required
just to prove that $\Ind(D_V)=\Ind(D_{\wV})$, and thus to deduce
the classical result that the arithmetic genus is birationally 
invariant.

\begin{thm}
\label{thm:holAbundle}
Let $M$ be a compact complex manifold, let $A$ be a unital $C^*$-algebra,
and let $\cV$ be an $A$-vector bundle which has a holomorphic structure.
Then $H^0(M, \cO_{\cV})$, the space of holomorphic sections of $\cV$,
is finitely generated and projective as an $A$-module.
\end{thm}
\begin{proof}
We know that the Dolbeault operator $D_M = \bar\partial +
{\bar\partial}^*$ is elliptic since $M$ is compact, and thus
$D_{\cV}$, or $D_M$
with coefficients in $\cV$ (defined with respect to some choice of
connection on $\cV$), is elliptic in the sense of
Mish\-chenko-Fo\-men\-ko. Thus the kernel of $D_{\cV}$,
and hence $H^0(M, \cO_{\cV})$ (which is the intersection of that
kernel with the $0$-forms with values in $\cV$) is a submodule of a
finitely generated projective $A$-module. It suffices to show that
$H^0(M, \cO_{\cV})$ is $A$-projective. This follows from the existence
of a ``reproducing kernel,'' which in turn follows from the Cauchy
integral formula (see \cite[Ch.\ 1, Proposition 2]{Nara}. 
(Locally, projection from the Hilbert $A$-module
of $L^2$-sections of $\cV$ to the holomorphic sections is given by the
same formula as projection from $L^2$-functions on a ball in $\bC^n$
to the Bergman space of
$L^2$-holomorphic functions, except that everything is now
$A$-linear.) Existence of the projection proves that $H^0(M,
\cO_{\cV})$  is a Hilbert $A$-module. Since it is also a submodule of
a finitely generated and projective $A$-module, it is itself finitely
generated and projective.
\end{proof}
\begin{thm}
\label{thm:KaehlAbundle}
Let $M$ be a compact K\"ahler manifold, let 
$\pi$ be a countable discrete group, let $u\co M\to B\pi$,
let $A=C^*(\pi)$ be a unital $C^*$-algebra,
and let $\cV=\cV_M$ be the associated flat $A$-bundle as above.
Then the $A$-index of the Dolbeault operator on $M$ with coefficients
in $\cV$ is equal to $\sum_{j=i}^n (-1)^j \left[ 
H^0(M, \Omega^j_{\cV})\right]$, the alternating sum of the
classes of the spaces of holomorphic $j$-forms with coefficients
in $\cV$, the sum being computed in $K_0(A)$. {\lp}The terms in the
sum are well-defined in $K_0(A)$ because of Theorem
\ref{thm:holAbundle}.{\rp} 
\end{thm}
\begin{proof}
By repeated application of Theorem \ref{thm:holAbundle}, each 
$H^0(M, \Omega^j_{\cV})$ is finitely generated and projective as an
$A$-module. On the other hand, complex conjugation exchanges
$\partial$ and $\bar\partial$ and preserves the bundle $\cV$
(which is just the complexification of the universal flat bundle for
the \emph{real} group {\Ca}). Because of the K\"ahler condition,
it sends $H^0(M, \Omega^j_{\cV})$ to the space of harmonic
$(0,j)$-forms with values in $\cV$ (see \cite[Theorem
4.7]{MR608414}, the proof of
\cite[Theorem 5.1]{MR608414}, and \cite[\S15.7]{Hir}), and thus this
latter space is also  
finitely generated and $A$-projective. In other words, compact
perturbation of the operation is not needed in defining the
Mish\-chenko-Fo\-men\-ko $A$-index of the Dolbeault operator on $M$
with coefficients in $\cV$, and thus the latter agrees with the
alternating sum $\sum_{j=i}^n (-1)^j \left[ 
H^0(M, \Omega^j_{\cV})\right]$.
\end{proof}

Theorems \ref{thm:holAbundle} and \ref{thm:KaehlAbundle} are the
analytical steps in the proof of Theorem
\ref{thm:indexinvariance}. To conclude the proof, there is one more
algebraic geometry or complex analysis part of the proof, namely:
\begin{thm}
\label{thm:blowingholforms}
Let $\phi\co \wV \to V$ be a blowing down of nonsingular 
$n$-dimensional complex
projective varieties as \textup{\S\ref{sec:proofmethod}} above, 
and let $\pi$ be a discrete group. Then for each $j$,
pull-back of forms induces an isomorphism
$H^0(V, \Omega^j_{\cV_V}) \xrightarrow{\cong}
H^0(\wV, \Omega^j_{\cV_{\wV}})$. 
\end{thm}
\begin{proof}
(Compare, for example, \cite[\S5.4]{MR637060} and
\cite[Proposition 6.16]{MR0453732}.)
Recall that $\phi$ is an isomorphism from $\wV\smallsetminus \wZ$
to $V\smallsetminus Z$. Furthermore, any holomorphic form on $V$
(resp., $\wV$) is determined by its restriction to the
dense subset $V\smallsetminus Z$ (resp., $\wV\smallsetminus \wZ$). So
$\phi^*$ is certainly injective. (If $\phi^*\omega = 0$,
then $\phi^*\omega |_{\wV\smallsetminus \wZ} = 0$, hence
$\omega|_{V\smallsetminus Z} = 0$, hence $\omega=0$.) So
we need only show that $\phi^*$ is surjective. Suppose $\eta$
is a holomorphic $j$-form on $\wV$ with values in $\cV_{\wV}$.
We need to show it comes from such a form on $V$. Since
$\phi$ is an isomorphism from $\wV\smallsetminus \wZ$
to $V\smallsetminus Z$ and $\phi$ induces an isomorphism
on $\pi_1$ (Proposition \ref{prop:IsoOnPi}), we may identify
$\eta_{\wV\smallsetminus \wZ}$ with a holomorphic $j$-form
on $V\smallsetminus Z$ with values in $\cV_V$. We just need
to show that it extends to all of $V$. (Then the pull-back of the
extension will be a holomorphic $j$-form that agrees with $\eta$
on a dense set, and hence is equal to $\eta$.) But $\eta$
is given to be holomorphic on the complement  in $V$ of the subvariety
$Z$, and it is locally bounded since it came from a continuous
form on the compact manifold $\wV$. So it extends by the
Riemann Extension Theorem \cite[Ch.\ 4, Proposition 2]{Nara}
or by related extension theorems such as Hartogs' Theorem
(see the rest of \cite[Ch.\ 4]{Nara}).
\end{proof}
Now Theorem \ref{thm:indexinvariance}, and thus Theorem
\ref{thm:main}, follows from Theorems \ref{thm:KaehlAbundle} and
\ref{thm:blowingholforms}.  

\section{A ``Converse Theorem''}
\label{sec:remarks}

Surgery theory can be used to prove that the signature is essentially
the only Pontrjagin number for closed oriented manifolds
(characteristic number of the form $\langle f(p_1,\,p_2,\,\cdots),
[M]\rangle$), with $f$ a polynomial
function of the rational Pontrjagin classes $p_1,\,p_2,\,\cdots$,
that is a homotopy invariant. In fact, the result can be improved for
non-simply connected manifolds, and the higher signatures are
essentially the only more general characteristic numbers of the
form $\langle f(p_1,\,p_2,\,\cdots)\cup u^*(x), [M]\rangle$
that can be homotopy invariant for all
$x\in H^*(B\pi_1(M),\bQ)$, i.e., are the only
characteristic functions of the form $u_*( f(p_1,\,p_2,\,\cdots)\cap [M])
\in H_*(B\pi_1(M),\bQ)$ that can be homotopy invariant.
A precise theorem to this effect may
be found, for instance, in \cite[Theorem 6.5]{Davis} or in
\cite[\S1.2, Remark 1.9]{KrL}. 
What Davis's theorem (\emph{loc.\ cit.}) really says is that
one has freedom to vary $\cL(M)$ more or less arbitrarily within
the oriented homotopy class of $M$, as long as one doesn't
change $u_*(\cL(M)\cap [M]) \in H_*(B\pi_1(M),\bQ)$.
But the rational Pontrjagin classes can be recovered from $\cL$,
so this implies that $p_1,\,p_2,\,\cdots$ can be varied
more or less arbitrarily within
the oriented homotopy class of $M$, as long as one doesn't
change $u_*(\cL(M)\cap [M]) \in H_*(B\pi_1(M),\bQ)$.  Another
version of a very similar result is the following:

\begin{thm}
\label{thm:biratPontnumbers}
If $\xi\co\Omega_k(B\pi)\otimes \bQ \to \bQ$ is a linear functional on the
rational oriented bordism of a group $\pi$, and if $\xi([M \xrightarrow{u}
B\pi])=\xi([\wM \xrightarrow{u\circ h} B\pi])$ whenever
$h\co\wM\to M$ is an orientation-preserving homotopy equivalence,
then $\xi$ is given by a higher signature, i.e.,
\[ \xi([M \xrightarrow{u}B\pi]) = \langle u_*(\cL(M)\cap [M]), x\rangle
\]
for some $x\in H^*(B\pi,\bQ)$.
\end{thm}
\begin{proof}
We know that 
\[ \Omega_k(B\pi)\otimes \bQ \cong \bigoplus_{0\le j\le \lfloor
\frac k4 \rfloor}
H_{k-4j}(B\pi,\bQ) \otimes_{\bQ} (\Omega_{4j}
\otimes_{\bQ}).
\]
Thus $\xi$ can be written as a sum $\sum_j x_j\otimes f_j$,
where $x_j\in H^{k-4j}(B\pi,\bQ) $ and $f_j$ is a linear functional on
$\Omega_{4j} \otimes {\bQ}$. The splitting of $\Omega_k(B\pi)\otimes
\bQ $ corresponds geometrically to a splitting of any bordism 
class into a sum of classes with representatives of the
form $M^{k-4j} \times N^{4j} \xrightarrow{u_j} B\pi$, where $u_j$
factors through $M^{k-4j}$ and $M^{k-4j} \xrightarrow{u_j} B\pi$
represents a class in $H_{k-4j}(B\pi,\bQ)$. So homotopy invariance of
$\xi$ means that each $f_j$ must be homotopy invariant on $\Omega_{4j}
\otimes {\bQ}$.  However, it is shown in \cite{MR0172306} that the
differences $[N] - [N']$, where $N$ and $N'$ are related by an
orientation-preserving homotopy equivalence, generate an ideal $I_*$ of
codimension $1$ in $\Omega_* \otimes {\bQ}$, with the quotient
$(\Omega_* \otimes {\bQ})/I_*$ detected by the signature. So $f_j$
must be a multiple of the signature. Absorbing the constant factor
into $x_j$, the Hirzebruch Signature Theorem proves that 
\[
\xi([M \xrightarrow{u}B\pi]) = \langle u_*(\cL(M)\cap [M]), x\rangle,
\quad \text{with } x = \sum x_j.
\]
\end{proof}

Because of Theorem \ref{thm:biratPontnumbers},
the Novikov Conjecture asserts that every functional on
$\Omega_*(B\pi)\otimes \bQ$ that \emph{could be} homotopy invariant
\emph{is} homotopy invariant.
The analogy which we proposed in Section \ref{sec:intro}, 
which involves replacing closed oriented manifolds by nonsingular
complex projective varieties, Pontrjagin classes by Chern classes, the
$L$-class by the Todd class, and homotopy equivalence by rational
equivalence, motivates the following analogues of this result:

\begin{thm}[cf.\ \cite{MR0172306}]
\label{thm:biratChernnumbers}
The differences $[N] - [N']$, where $N$ and $N'$ are
birationally equivalent nonsingular complex projective varieties,
generate an ideal $I_*$ of codimension $1$ in the complex 
cobordism ring  $\Omega_*^U \otimes {\bQ}$, with the quotient
$(\Omega_*^U \otimes {\bQ})/I_*$ detected by the Todd genus.
\end{thm}
\begin{proof}
Clearly these differences generate an ideal, since 
if $\xymatrix{N\ar@{.>}^h[r] & N'}$ is a birational map, so is 
$\xymatrix{N\times P\ar@{.>}^{h\times \text{id}}[r] & N'\times P}$,
for any nonsingular complex projective variety $P$, and also all of
$\Omega_*^U \otimes  {\bQ}$ ($*>0$) is represented by nonsingular complex
projective varieties \cite[p.\ 130]{Stong}.
In fact one can take the complex projective spaces as
generators for $\Omega_*^U
\otimes {\bQ}$. It's non-trivial that one has inverses within the
monoid generated by the nonsingular projective varieties, but it's
not known if one can choose the varieties to be connected. Fortunately this
last point doesn't concern us.  In any event, it suffices to observe
that we have a $\bQ$-basis for $\Omega_{2k}^U\otimes \bQ$ consisting
of all products $\bC\bP^{j_1}\times \bC\bP^{j_2}\times \cdots \times
\bC\bP^{j_r}$ 
with $j_1\ge j_2\ge \cdots \ge j_r$ and $j_1 + j_2 + \cdots + j_r = k$
\cite[page 111]{Stong}.
But all of these varieties are in the same birational equivalence class,
since they all have the same function field, namely a purely transcendental
extension $\bC(x_1,\cdots,x_k)$ of $\bC$ of transcendence degree $k$.
(To put things another way, all the above products of projective
spaces are clearly projective completions of the same affine $k$-space.)
Thus $I_{2k}$ has codimension $1$ in $\Omega_{2k}^U\otimes \bQ$,
and the quotient is detected by the Todd genus, since all of these
products of projective spaces have Todd genus $1$. That completes the proof.
\end{proof}
\begin{thm}
\label{thm:birathigherChernnumbers}
Suppose $\pi$ is a group such that for all $k$, $H_{2k}(B\pi,\bQ)$ is 
spanned by the classes of maps from $k$-dimensional
nonsingular projective varieties into $B\pi$.
{\lp}Just as an example, this is certainly the case if $\pi$
is free abelian, since then all even homology of $\pi$ is
generated by $k$-dimensional abelian varieties.{\rp}
If $\xi\co\Omega_{2k}^U(B\pi)\otimes \bQ \to \bQ$ is a linear functional on the
rational complex bordism of a group $\pi$, and if $\xi([M \xrightarrow{u}
B\pi])=\xi([\wM \xrightarrow{u\circ h} B\pi])$ whenever
$h\co\wM\to M$ is a{\lp}n everywhere defined{\rp} birational
equivalence of nonsingular complex projective varieties
of complex dimension $k$, 
then $\xi$ is given by a higher Todd genus, i.e.,
\[ \xi([M \xrightarrow{u}B\pi]) = \langle u_*(\cT(M)\cap [M]), x\rangle
\]
for some $x\in H^*(B\pi,\bQ)$.
\end{thm}
\begin{proof}
We can prove this the same way we proved Theorem \ref{thm:biratPontnumbers},
with Theorem  \ref{thm:biratChernnumbers} substituting for Kahn's
Theorem \cite{MR0172306}.  Only one point requires clarification:
the extra hypothesis on $\pi$ that $H_{2k}(B\pi,\bQ)$ is
spanned by the classes of maps from $k$-dimensional
nonsingular projective varieties into $B\pi$. This ensures that
$\Omega_{2k}^U(B\pi)\otimes \bQ$ is spanned by classes of
maps $ M^{2k-2j} \times N^{2j} \xrightarrow{u_j} B\pi$, where $u_j$
factors through $M^{2k-2j}$, $M^{2k-2j} \xrightarrow{u_j} B\pi$
represents a class in $H_{2k-2j}(B\pi,\bQ)$, and $M$ and $N$ are
\emph{smooth projective varieties} of complex dimensions $k-j$ and $j$,
respectively. (Otherwise we would just know that we could choose
$M$ to be stably almost complex, which isn't good enough for our
purposes.) The rest of the proof is the same as before.
\end{proof}

\section{Examples and Concluding Remarks}
\label{sec:concl}

To conclude, we mention a few examples where the hypotheses of
Theorems \ref{thm:main} and \ref{thm:birathigherChernnumbers} apply, and
explicate exactly what these theorems say in these special cases.
We also mention some open problems.

\begin{example}
The simplest case of our theory is when $\pi$ is free abelian. Since
the first Betti number of a compact K\"ahler manifold has to
be even, it is natural to start with the case $\pi=\bZ^2$, which of
course is the fundamental group of an elliptic curve $E$. The
hypotheses of Theorems \ref{thm:main} and
\ref{thm:birathigherChernnumbers} certainly apply; in fact,
essentially all known ways of proving SNC apply in this case. Putting
Theorems \ref{thm:main} and \ref{thm:birathigherChernnumbers}
together, we conclude the following:
\begin{enumerate}
\item Let $V$ be a smooth projective variety of complex dimension $n
  \ge 1$ with a homomorphism $u\co \pi_1(V)\to \bZ^2$. Then $\langle
  \cT(V)\cup u^*(x),\,[V]\rangle$ is a birational invariant of $V$,
  for $x$ the usual generator of $H^2(B\bZ^2,\bZ) \cong \bZ$. (This is
  the only non-trivial higher Todd genus in this case.)
\item \emph{Every} birationally invariant linear functional on
  $\Omega_{2n}^U (B\bZ^2) \otimes \bQ$ is a linear combination of the
  Todd genus (or arithmetic genus) and the higher Todd genus of (1).
\end{enumerate}
For example, suppose we consider the case of surfaces ($n=2$). In this
case, the Todd genus is $\langle \frac{1}{12}(c_1(V)^2 + c_2(V)),\,
[V]\rangle$, but we see that $\langle c_1\cup  u^*(x),\,[V]\rangle$ is
also a birational invariant. In particular, suppose $V$ is birationally
equivalent to $E \times \bP^1$ (in which case we can choose $u$ to be an
isomorphism on $\pi_1$). Now $H^*(E \times \bP^1,\bZ)$ is the
tensor product of an exterior algebra on generators $x_1$ and $x_2$
(each of degree $1$) and a truncated polynomial ring on a generator
$y$ in degree $2$, and we can choose $u^*(x) = x_1\cup x_2$, and the
Chern classes of $E \times \bP^1$ are $c_1 = 2y$, $c_2 = 0$. So $V$
has vanishing Todd genus, hence $c_1(V)^2 + c_2(V) = 0$, and
\[
\langle c_1(V)\cup  u^*(x),\,[V]\rangle =
\langle c_1(E \times \bP^1)\cup  u^*(x),\,[E \times \bP^1]\rangle 
= \langle x_1\cup x_2 \cup 2y,\,[E \times \bP^1]\rangle = 2.
\]
This implies (among other things) that $c_1(V)$ has to be non-zero and
even, which we wouldn't know just from the vanishing of the Todd
genus. For example, consider the case where $V$ is $E \times \bP^1$
with one point blown up. Topologically, such a $V$ is 
a connected sum $(T^2\times S^2) \# \overline{\bC\bP^2}$, which has
Euler characteristic $1$ and signature $-1$. So $\langle c_2(V),\,
[V]\rangle = 1$ (the top Chern class is also the Euler class) and
$\langle c_1(V)^2,\,[V]\rangle = -1$, which we could deduce in many
ways, but more interestingly, $\langle c_1(V)\cup x_1\cup x_2
,\,[V]\rangle =  2$. (In fact, since the cohomology ring of $V$ is
the same as for $E \times \bP^1$, but with one additional generator
$z$ in degree $2$ with $z^2 = - x_1\cup x_2 \cup y$, $z\cup x_1 = 0$,
$z\cup x_2 = 0$, $z\cup y = 0$, we find that
$c_2(V) = - z^2$, $c_1(V) = 2y + z$, which does satisfy the above
conditions.)
\end{example}

\begin{example}
For another example, suppose $V_1$ and $V_2$ are birationally
equivalent and both of them are aspherical, i.e., they are each
homotopy equivalent to $B\pi$, where $\pi$ is the common fundamental
group. We fix the isomorphism $\pi_1(V_1) \to \pi_1(V_2) \to \pi$ determined
by the birational equivalence; this fixes an identification of both
cohomology rings $H^*(V_1,\bZ)$ and $H^*(V_2,\bZ)$ with the group
cohomology ring $H^*(\pi,\bZ)$. We also assume that SNC holds for
$\pi$; this is the case, for instance, if either $V_1$ or $V_2$ is
locally symmetric, since then $\pi$ is a discrete subgroup of a Lie
group and \cite{K1} applies. (In fact, \cite{KS} also applies in this
case.) The Novikov Conjecture implies that $\cL(V_1) = \cL(V_2)$, and
thus $V_1$ and $V_2$ have the same rational Pontrjagin classes. But
what about the Chern classes? Some constraints on them are implied by
the relationship between the Chern and Pontrjagin classes, that
if (formally) $1+c_1+c_2+\cdots = \prod (1+x_i)$, then
$1+p_1+p_2+\cdots = \prod (1+x_i^2)$. But this by itself does
not determine the Chern classes.  (For example, we have $p_1 = c_1^2 -
2c_2$, and in the case of an algebraic surface, $c_2$ is
determined by the Euler characteristic, so we know $c_1^2$, but not
necessarily $c_1$ itself.) Conjecture \ref{conj:higherTodd}, which
in this case is an actual theorem, thanks to 
Theorem \ref{thm:main}, says
that we also have $\cT(V_1) = \cT(V_2)$. Thus, for example,
$c_1(V_1) = c_1(V_2)$, which is an extra piece of information. I am
not sure if it is possible to deduce that all of the Chern classes
of $V_1$ and $V_2$
must agree, but this seems plausible, and it certainly holds if
$\dim_\bC V \le 3$. For example, if $\dim_\bC V_1 = \dim_\bC V_2 = 3$,
then $c_3(V_1)=c_3(V_2)$ since the Euler characteristics agree, and as
we have seen, $c_1(V_1)=c_1(V_2)$ since the part of $\cT$ in degree
$2$ must match up, and then $c_2(V_1)=c_2(V_2)$ since $2c_2 = c_1^2 -
p_1$. 
\end{example}

\begin{rem}
An obvious open problem is to find a purely algebraic analogue of
these results, preferably one which would also apply in characteristic
$p$. This may require a new notion of ``fundamental group,'' since the
usual algebraic fundamental group only takes finite coverings into
account, and if $\pi$ is finite, its rational cohomology vanishes, so
that the higher Todd genera do not tell us any more than the Todd
genus itself.
\end{rem}

\renewcommand{\MR}{\relax}
\bibliography{higherTodd}
\bibliographystyle{amsplain}

 \end{document}